\documentclass[12pt]{article}
\usepackage{amsmath,amssymb}
\usepackage{amsthm}
\usepackage{cite}

\newtheorem{theorem}{Theorem}

\newtheorem{lemma}{Lemma}

\newtheorem{prop}{Proposition}

\def\PC{\mathcal{P}}
\def\MC{\mathcal{M}}

\def\LC{\mathcal{L}}

\def\B{\mathbf{B}}

\def\E{\mathbf{E}}

\def\R{\mathbf{R}}

\def\1{\mathbf{1}}

\def\al{\alpha}
\def\be{\beta}
\def\pa{\partial}
\def\ep{\epsilon}
\def\de{\delta}

\def\ka{\varkappa}

\newcommand{\la}{\lambda}
\newcommand{\si}{\sigma}

\newcommand{\Om}{\Omega}

\begin{document}
	
\title{Regularity and Sensitivity for McKean-Vlasov SPDEs
}
\author{
	Vassili N. Kolokoltsov\thanks{University of Warwick,
		Coventry CV4 7AL UK, IPI RAN RF and St.-Petersburg State Univ., Email: v.kolokoltsov@warwick.ac.uk}
	and Marianna Troeva\thanks{Research Institute of Mathematics, North-Eastern Federal University,
		58 Belinskogo str., Yakutsk 677000 Russia, Email: troeva@mail.ru
	}}
	\maketitle

	\maketitle

	\begin{abstract}
		In the first part of the paper we develop first the sensitivity analysis for the nonlinear McKean-Vlasov diffusions stressing precise estimates of growth of the solutions and their derivatives with respect to the initial data, under	rather general assumptions on the coefficients. The exact estimates become particularly important when treating the extension of these equations having random coefficient, since the noise is usually assumed to be unbounded.
		The second part contains our main results dealing with the sensitivity of stochastic  McKean-Vlasov diffusions. By using the method of stochastic characteristics, we transfer these equations to the non-stochastic equations with random coefficients thus making it possible to use the estimates obtained in the first part. The motivation for studying sensitivity of McKean-Vlasov SDEs arises naturally from the analysis of the mean-field games with common noise.
	\end{abstract}

	

\medskip


	\section{Introduction}

{\it Nonlinear diffusion equations}\index{nonlinear diffusion}
 represent the general class of equations of diffusion type with coefficients depending on the unknown function.
 In other words, they are evolutionary equations of the second order, which are linear with respect to the derivatives:
\begin{equation}
\label{eqdefgenMcKeanVl}
\frac{\pa u}{\pa t} =\frac12 (A(u,x) \nabla,\nabla)  u(x)+(b(u,x),\nabla ) u(x)+V(u,x) u,
\end{equation}
or more explicitly
\begin{equation}
\label{eqdefgenMcKeanVl1}
\frac{\pa u}{\pa t} =\frac12 \sum_{i,j}A_{ij}(u,x) \frac{\pa ^2 u(x)}{\pa x_i \pa x_j}
+ \sum_j b_j(u,x) \frac{\pa u(x)}{\pa x_j} +V(u,x) u(x),
\end{equation}
with given functions $A(u,x)$, $b(u,x)$, $V(u,x)$.

Very often these nonlinear equations arise in their weak form, as equations of type
\begin{equation}
\label{eqgennonlinearprobrep}
 \frac{d}{dt}(f,\mu_t)=(L(\mu_t)f,\mu_t), \quad \mu_0=Y, \quad f\in D,
\end{equation}
where
\begin{equation}
\label{eqdifnonlgen}
L(\mu)f(x) =\frac12  (A(\mu,x)\nabla,\nabla)f(x) + (b(\mu, x),\nabla) f(x)+V(\mu, x) f(x).
\end{equation}

The strong form of equation (\ref{eqgennonlinearprobrep}) written for measures $\mu$ having densities $\phi$ is
\begin{equation}
\label{eqdifnonlgen1}
\frac{\pa \phi}{\pa t} =\frac12 ( \nabla, \nabla (A(\phi,x) \phi (x))) - (\nabla, (b(\phi, x)\phi(x)))+V(\phi, x) \phi(x),
\end{equation}
where we denote the functions of measures and of their densities by the same letters, with some abuse of notations. More explicitly,
equation (\ref{eqdifnonlgen1}) rewrites as
\begin{equation}
\label{eqdifnonlgen2}
\frac{\pa \phi}{\pa t} =\frac12  \sum_{i,j}\frac{\pa ^2}{\pa x_i \pa x_j}(A_{ij}(\phi,x)\phi)
- \sum_j \frac{\pa }{\pa x_j} (b(\phi, x))\phi(x))+V(\phi, x) \phi(x),
\end{equation}
which can be of course again rewritten in general form (\ref{eqdefgenMcKeanVl1}) (but with some other functions $A$, $b$ and $V$).

Nonlinear diffusion equations are often referred to as the {\it McKean-Vlasov diffusions}\index{McKean-Vlasov diffusions},
especially when the coefficients $A,b,V$ depends on the unknown function $u$ via the moments of the type
$\int g(x_1, \cdots, x_k)u(x_1) \cdots u(x_k) dx_1 \cdots dx_k$.

In case of the matrix $A(u,x)$ not depending on $u$, equation (\ref{eqdefgenMcKeanVl})
belongs to the class of equations of the type
\begin{equation}
\label{eqdefgenMcKeanVl1a}
\frac{\pa u}{\pa t} =A u(x)+(b(u,x),\nabla u(x))+V(u,x) u,
\end{equation}
where $A$ is a usual second order operator, which we shall treat here.
More precisely, we shall work with the corresponding weak equation (\ref{eqgennonlinearprobrep}) and the corresponding strong form
\begin{equation}
\label{eqdifnonlgen11}
\frac{\pa \phi}{\pa t} =\frac12 ( \nabla, \nabla (A(x) \phi (x))) - (\nabla, (b(\phi, x)\phi(x)))+V(\phi, x) \phi(x).
\end{equation}

The stochastic  McKean-Vlasov equations have the form
\begin{equation}
\label{eqlimMacVlaSPDE}
d(f, \mu_t)=(L_{t}(\mu_t)f, \mu_t) \, dt +(\si_{com}(.)\nabla f, \mu_t) \, dW_t,
\end{equation}
which is written here in the weak form meaning that it should hold for all $f \in C^2(\R^d)$.
Here
\begin{equation}
\label{eqlimMacVlaSPDE1}
L_{t}(\mu_t)f (x)= \frac12
((\sigma_{com}\sigma_{com}^T+\sigma_{ind}\sigma_{ind}^T)(x)\nabla,\nabla)f (x)
+(b(t,x, \mu_t),\nabla)f (x).
\end{equation}

The motivation for studying sensitivity of McKean-Vlasov SDEs (\ref{eqlimMacVlaSPDE}) and the notations $\sigma_{com}, \sigma_{ind}$
arise naturally from the analysis of the mean-field games with common noise in which positions of $N$ agents are governed by the system of SDEs
\begin{equation}
	\label{eqstartSDEN}
	dX_t^i=b(t,X_t^i, \mu_t^N, u_t^i)\, dt +\si_{ind}(X_t^i) dB_t^i +
	\si_{com}(X_t^i) dW_t,
\end{equation}
where all $X_t^i$ belong to $\R^d$, $B_t^1, \cdots, B_t^N$ are independent $d'$-dimensional standard Brownian motions, $W_t$ is a one-dimensional standard Brownian motion. $W_t$
referred to as the common noise, and all $B^j_t$, referred to as the idiosyncratic or individual noises. The subscripts 'com' and 'ind' referred to the objects related to the common or to the individual noises.
The parameters $u_t^i\in U\subset \R^m$ are controls available to the players,
trying to minimize their payoffs.

Let us mention directly that in our approach it is more convenient to work with the Stratonovich
differentials. Namely, by the usual rule $Y\circ dX=Y dX +\frac12 dY dX$, equation (\ref{eqlimMacVlaSPDE})
rewrites in the Stratonovich form as
\begin{equation}
\label{eqlimMacVlaSPDEst}
d(f, \mu_t)=(L_t^{St}(\mu_t) f, \mu_t) \, dt +
(\si_{com}(.) \nabla f^T , \mu_t) \, \circ dW_t
\end{equation}
where
\[
L_t^{St}(\mu_t)f (x)=
\frac{1}{2}(\si_{ind}(x)\si_{ind}^T(x)\nabla,\nabla)f(x)
+\left((b(t,x,\mu_t)-\frac{1}{2}(\si_{com}(x)\nabla \si_{com}(x))),\nabla \right)f(x),
\]
where we denote $(\si_{com}(x)\nabla \si_{com}(x))=\left((\si_{com}(x),\,\nabla) \si^1_{com}(x),\dots,(\si_{com}(x),\,\nabla) \si^d_{com}(x)\right)$.

It is known (see e.g. \cite{KurXi}) that, for fixed common functions
$u^i_t(X_t^i)=u_t(X_t^i)$, and under appropriate regularity assumptions on $b,\si_{ind},\si_{com}$
the system (\ref{eqstartSDEN}) is well-posed and the corresponding empirical measures $\mu_t^N$
converge, as $N\to \infty$, to the unique solution $\mu_t$ of the nonlinear SPDE of the McKean-Vlasov type.

Notice that there is quite an extensive literature on properties of McKean-Vlasov SPDEs  (see e.g. \cite{KurXi, DaVa95, HuNu13, CoJoKh13, KolTro15, CrMcMur} and references therein). The well-posedness of the McKean-Vlasov SPDE was shown in \cite{KurXi} in the class of $L_2$-functions, and for measures in \cite{DaVa95}, though under an additional monotonicity assumption.

Our paper is organized as follows.
In the first part of the paper we develop first the sensitivity analysis for the nonlinear McKean-Vlasov diffusions stressing precise estimates of growth of the solutions and their derivatives with respect to the initial data, under	rather general assumptions on the coefficients. The exact estimates become particularly important when treating the extension of these equations having random coefficient, since the noise is usually assumed to be unbounded.
	
The second part contains our main results dealing with the sensitivity of stochastic  McKean-Vlasov diffusions.
Our basic approach will be the method of stochastic characteristics, see \cite{Kunita}, \cite{KoTyu03}, though in its simplest form, available for either one-dimensional noise or a multidimensional white noise with constant correlations. This method allows one to turn stochastic  McKean-Vlasov equation into a non-stochastic equation of the second order, but with random coefficients, and thus to use the estimates obtained in the first part.

We shall be a bit sketchy in the development of the first part, as the results here are already essentially known, and only stress the way the explicit bounds are obtained.

The following basic notations will be used:

$C^n(\mathbf{R}^d)$ is a Banach space of $n$ times continuously differentiable and bounded functions $f$ on $\mathbf{R}^d$  such that each derivative up to and including order $n$ is bounded, equipped with norm $\|f\|_{C^n}$ which is the supremum of the sums of all the derivatives up to and including order $n$.

$C_{\infty}(\R^d)$ is a Banach space of bounded continuous functions $f:\R^d\to \R$ with $ \lim_{x\to \infty}f(x)=0$, equipped with sup-norm.

$C_{\infty}^n(\mathbf{R}^d)$ is a closed subspace of $C^n(\mathbf{R}^d)$ with $f$ and all its derivatives up to and including order $n$ belonging to $C_{\infty}(\mathbf{R}^d)$.

If $\MC$ is a closed subset of a Banach space $\B$, then

$C([0,T], \MC)$ is a metric space of continuous functions $t \rightarrow \mu_t \in \MC$
with distance $\|\eta-\xi\|_{C([0,T], \MC)}=\sup_{t\in [0,T]} \|\eta_t-\xi_t\|_\B$. An element from $C([0,T], \MC)$ is written as $\{\mu_.\}=\{\mu_t, t\in[0,T]\}$.


$\MC (\R^d)$ is a Banach space of finite signed Borel measures on  $\R^d$.

$\MC^+(\R^d)$ and  $\PC (\R^d)$ the subsets of $\MC (\R^d)$ of positive and positive normalised (probability) measures, respectively.

Let $\MC_{<\la}(\R^d)$ (resp. $\MC_{\le \la}(\R^d)$
 or $\MC_{\la}(\R^d)$) and $\MC_{<\la}^+(\R^d)$ (resp. $\MC_{\le \la}^+(\R^d)$ or $\MC_{\la}(\R^d)$)
denote the parts of these sets
containing measures of the norm less than $\la$ (resp. not exceeding $\la$ or equal $\la$).

Let $C^{k\times k}(\R^{2d})$ denote the subspace of $C(\R^{2d})$ consisting of functions $f$ such that the partial derivatives $\pa^{\al+\be}f/\pa x^{\al} \pa y^{\be}$ with multi-index $\al,\be$, $|\al|\le k, |\be|\le k$, are well defined and belong to $C(\R^{2d})$. Supremum of the norms of these derivatives provide the natural norm for this space.

%

For a function $F$ on $\MC^+_{\le \la}(\R^d)$ or $\MC_{\le \la}(\R^d)$ the variational derivative is defined as the directional derivative of $F(\mu)$ in the direction $\delta_x$:
\[
\frac{\de F(\mu)}{\de \mu(x)}=\frac{d}{dh}|_{h=0}F(\mu +h\de x).
\]
The higher derivatives
$\delta^l F(\mu)/\de \mu(x_1)...\de \mu(x_l)$ are defined inductively.


Let $C^k(\MC_{\le \la}(\R^d))$ denote the space of functionals such that the $k$th order variational derivatives are well defined
and represent continuous functions of all variables with measures considered in their weak topology

Let $C^{k,l}(\MC_{\le \la}(\R^d))$ denote the subspace of $C^k(\MC_{\la}(\R^d))$ such that all derivatives up to order
$k$ have continuous bounded derivatives up to order $l$ as functions of their spatial variables.

let $C^{2,k\times k}(\MC_{\le \la}(\R^d))$ be the space of functionals such that their second order variational derivatives
are continuous as functions of all variables and belong to
$C^{k\times k}(\R^{2d})$ as functions of the spatial variable; the norm of this space is
\[
\|F\|_{C^{2,k\times k}(\MC_{\le \la}(\R^d))}
=\sup_{\mu\in \MC_{\la}(\R^d)}\left \|\frac{\de^2 F}{\de \mu (.) \de \mu(.)}\right\|_{C^{k \times k}(\R^{2d})}.
\]

$(f,\mu)=\int f(x)\mu(dx)$ denotes the usual pairing of functions and measures on $\R^d$.

$\E$ denotes the expectation.

$E_{\beta}$ denotes the Mittag-Leffler function.

For mappings $F: M \to B_2$ from a closed convex subset $M$ of a Banach space $B_1$ to a Banach space $B_2$,
we define the spaces $C^1(M,B_2)$ and $C^2(M,B_2)$ of differentiable functions
$F$ of order $1$ or $2$ with continuous bounded derivatives, the continuity being understood as the continuity
 of mappings $Y\mapsto DF(Y)$ and $Y\mapsto D^2F(Y)$ with the norm topologies of $M$, $\LC(B_1,B_2)$ and $\LC^2(B_1,B_2)$.
The norm on these spaces is defined as follows:
\[
\| F\|_{C^2(M,B_2)}=\sup_Y \left( \|F(Y)\|_{B_2}+ \|DF(Y)\|_{\LC(B_1,B_2)}
+ \|D^2F(Y)\|_{\LC^2(B_1,B_2)}\right)
\]
\begin{equation}
\label{eqdefderbetweenBanach}
=\sup_Y \left( \|F(Y)\|_{B_2}+ \sup_{\|\xi\|_{B_1}\le 1} \|DF(Y)[\xi]\|_{B_2}
+ \sup_{\|\xi\|_{B_1}, \|\eta\|_{B_1}\le 1} \|D^2F(Y)[\xi, \eta]\|_{B_2}\right).
\end{equation}

Weakening the condition of Lipschitz continuity it is handy to define the space $C^1_{luc}(M,B_2)$
(resp. $C^1_{uc}(M,B_2)$) of the mappings with locally uniformly continuous derivatives
(resp. uniformly continuous), which is the closed subspace $C^1(M,B_2)$ of functions $F$ such
that the mapping $Y\mapsto DF(Y)$ is uniformly continuous on bounded subsets of $M$ (resp. on the whole $M$).
If $F\in C^1_{luc}(M,B_2)$ (resp. $F\in C^1_{uc}(M,B_2)$), then 
\begin{equation}
\label{eqdefdirder3uni}
\|F(Y+\xi)-F(Y)-DF(Y)[\xi]\|_{B_2}\le \ep \|\xi\| \quad \mbox{for} \quad \| \xi \| \le \de,
\end{equation}
uniformly for $Y$ and $Y+\xi$ from any bounded set (resp. for any $Y,Y+\xi$ from $M$).

\section{Sensitivity for McKean-Vlasov Equations}

\begin{theorem}
\label{thwellposMcVla}
Let $b(\mu,.)=\{b^j(\mu,.)\}, V(\mu,.)$ be measurable bounded real functions with $V$
non-positive, satisfying (\ref{eq0thwellposfirstordernonl2}), (\ref{eq1thwellposfirstordernonl2})
\begin{equation}
\label{eq0thwellposfirstordernonl2}
V=\sup_{y,\mu} |V(\mu,y)| <\infty, \quad b=\sup_{y,\mu} \sum_j|b^j(\mu,y)| <\infty,
\end{equation}
and
\begin{equation}
\label{eq1thwellposfirstordernonl2}
\sup_x|b_t(\mu,x)-b_t(\eta,x)| \le L_A \|\mu-\eta\|_{\MC(\R^d)}, \quad \sup_x|V_t(\mu,x)-V_t(\eta,x)| \le L_A \|\mu-\eta\|_{\MC(\R^d)}.
\end{equation}
Let $A(x)$ be a matrix-valued function with elements belonging to $C^2(\R^d)$, which is uniformly elliptic, so that
\begin{equation}
\label{eq1thwellposMcVla}
m^{-1}(\xi,\xi) \le (A(x)\xi,\xi)\le m(\xi,\xi), \quad \sum_{ij}\|A_{ij}(x)\|_{C^2(\R^d)} \le M <\infty,
\end{equation}
with some positive constants $m,M$.

Then for any $T>0$ and any $Y$ having a nonnegative density $\phi$ (resp. without density), the mild
equation
\[
\phi_t(x)=\int G_t(x,y) Y(dy) + \int_0^t ds \int G_{(t-s)} (x,y) V_s(\phi_s , y)\phi_s(y) dy
\]
\begin{equation}
\label{eq1athwellposMcVla}
- \int_0^t ds \int \left(\frac{\pa}{\pa y} G_{(t-s)} (x,y), b_s(\phi_s , y) \phi_s (y) \right) dy
\end{equation}
has the unique bounded nonnegative solution
$\phi_. \in C([0,T], L_1(\R^d))$ (resp. unique solution $\phi_. \in C((0,T], L_1(\R^d))$ such that
 $\phi_t \to Y$ weakly, as $t \to 0$) such that
\[
\|\phi_t\|_{L_1(\R^d)}\le \|Y\|_{\MC(\R^d)},
\]
and for any two solutions $\phi_t^1$ and $\phi_t^2$
with the initial conditions $Y^1,Y^2$ the estimate
\begin{equation}
\label{eq2thwellposMcVla}
\|\phi_t^1 -\phi_t^2\|_{L_1(R^d)} \le \|Y^1-Y^2\|_{\MC(R^d)}E_{1/2} (\ka(T) \sqrt t)
\end{equation}
holds, where
\[
\ka(T) =  \bar C[(V\sqrt T+b)+L_A(\sqrt T+1)\|Y\|,
\]
with a constant $\bar C$ depending on $m,M,T$.

Finally, solutions $\phi_t$ to the mild equation (\ref{eq1athwellposMcVla}) also solve the corresponding weak equation
\begin{equation}
\label{eq4thwellposfirstordernonl}
\frac{d}{dt}(f,\mu_t)=\left(\frac12 (A \nabla, \nabla) f + (\nabla, b_t(\mu_t, x)f(x))+V_t(\mu_t, x)f(x), \mu_t\right),
\quad \mu_0=Y\in \MC^+_{\le \la}(\R^d).
\end{equation}
\end{theorem}

\noindent{\it Proof.} 
The solution $\phi_t$ is the fixed point of the mapping $\Phi_Y$, 
that is,
\[
\Phi_Y(\phi_.)(t)(x)=\int G_t(x,y) Y(dy) + \int_0^t ds \int G_{(t-s)} (x,y) V_s(\phi_s , y)\phi_s(y) dy
\]
\begin{equation}
\label{eq7thwellposfirstordernonlgen}
- \int_0^t ds \int \left(\frac{\pa}{\pa y} G_{(t-s)} (x,y), b_s(\phi_s , y) \phi_s (y) \right) dy.
\end{equation}

As is known from the theory of ordinary diffusions (see e.g. \cite{Porper, Ito}), there exist constants $\si$ and $C$ depending only on $m,M,T$ such that
the Green function  $G_t(x,y)$ of the operator $\phi(x) \mapsto \frac12 ( \nabla, \nabla (A(x) \phi (x)))$
(recall for clarity that it is just the transpose kernel to the Green function $G_t(y,x)$ for the dual operator
 $\frac12 (A(x)\nabla, \nabla)$) is differentiable in $x$ and $y$ and satisfies the estimate
\begin{equation}
\label{eq3thwellposMcVla}
 G_{t}(x,y) \le C G_{\si t}( x-y), \quad 0<t<T,
\end{equation}
\begin{equation}
\label{eqdiffeqheatkernelderiv}
\max \left(\left|\frac{\pa}{\pa y} G_{t}(x,y)\right|, \left|\frac{\pa}{\pa x} G_{t}(x, y)\right|\right)
\le Ct^{-1/2} G_{\si t}( x-y), \quad 0<t<T.
\end{equation}

Therefore, we can estimate the action of $\Phi^n$ as
\[
\|\Phi_Y(\phi_.)(t)\|_{L_1(\R^d)} \le \|Y\|_{L_1(\R^d)}+C(m,M,T) (V+b) \int_0^t   (t-s)^{-1/2} \|\phi_s\|_{L_1(\R^d)}\, ds.
\]
Iterating and using the definition of the Mittag-Leffler function yields
\begin{equation}
\label{eq3thwellposfirstordernonl2}
\|\Phi^n_Y(Y)(t)\|_{L_1(\R^d)} \le \|Y\|_{L_1(\R^d)} E_{1/2} (C(m,M,T)(V+b) \sqrt t),
\end{equation}
implying the uniform boundedness of all iterations.
Moreover,
\[
\|\Phi_Y(\phi_.^1)(t)-\Phi_Y(\phi_.^2)(t)\|_{L_1(\R^d)}
\le \ka(m,M,T) \int_0^t (t-s)^{-1/2} \|\phi_.^1 -\phi_.^2\|_{C([0,s],L_1(R^d))} \, ds,
\]
so that the existence and uniqueness of the fixed point follows from the basic fixed point principle in Banach spaces.

The fixed point principle can be applied to find the unique solution in the Banach space of measures $\MC(\R^d)$.
However, the curve $\Phi_Y(\mu_.)(t)$ is only weakly continuous in $t$ at $t=0$, because so is the first term in
(\ref{eq7thwellposfirstordernonlgen}). It follows from (\ref{eq7thwellposfirstordernonlgen}) that $\Phi_Y(\mu_.)(t)$
has a density for all $t>0$.

If  $\phi_t$ solves (\ref{eq7thwellposfirstordernonlgen}), then it solves the equation
\begin{equation}
\label{eqexermildweaklink}
(f, \phi_t) =(T_tf, \phi_t)+\int_0^t ds \left( [V_s(\phi_s,.) +(b_s(\phi_s,.),\nabla)] T_{t-s}f, \phi_t\right)
\end{equation}
for any $f\in C(\R^d)$,  where $T_t$ is the semigroup generated by the operator $(A \nabla, \nabla)/2$.
If $f\in C^2_{\infty}(\R^d)$, then we can differentiate this equation yielding (\ref{eq4thwellposfirstordernonl}) as required.
\hfill $\Box$ 

We can now obtain our first result on sensitivity.
\begin{theorem}
\label{thwellposMcVlasense}
Under the assumptions of Theorem \ref{thwellposMcVla} let the variational derivatives of $V_t(\mu,y)$
 and $b_t(\mu,y)$ with respect $\mu$ are well defined, are locally bounded and such that
 \begin{equation}
\label{eq0thwellposnonldif1}
\sup_{y,z,s} \left|\frac{\de V_s(\mu,y)}{\de \mu(z)} \right| \le R(\la ),
\quad \sup_{y,z,s} \sum_j\left|\frac{\de b^j_s(\mu,y)}{\de \mu(z)} \right| \le R(\la )
\end{equation}
for $\mu \in \MC^+_{\le \la}(\R^d)$ and constants $R(\la)$,
\begin{equation}
\label{eq1thwellposnonldif1}
\sup_{y,z,s} \left|\frac{\de V_s(\mu^1,y)}{\de \mu(z)} -\frac{\de V_s(\mu^2,y)}{\de \mu(z)} \right|
 \le L(\la)\|\mu^1-\mu^2\|_{L_1(\R^d)},
\end{equation}
\begin{equation}
\label{eq2thwellposnonldif1}
\sup_{y,z,s} \sum_j\left|\frac{\de b^j_s(\mu^1,y)}{\de \mu(z)} -\frac{\de b^j_s(\mu^2,y)}{\de \mu(z)} \right|
 \le L(\la)\|\mu^1-\mu^2\|_{L_1(\R^d)}
\end{equation}
for $\mu \in \MC_{\le \la}(\R^d)$ and constants $L(\la)$.

Then the mapping $\phi=\phi_0 \mapsto \phi_t$, with the initial data $\phi \in L_1(\R^d)$, or more generally $Y \mapsto \phi_t$
 with the initial data $Y \in \MC(\R^d)$, solving (\ref{eq7thwellposfirstordernonlgen}) according to Theorem \ref{thwellposMcVla},
belongs to $C^1_{luc}(L_1(\R^d),L_1(\R^d))$, or more generally to $C^1_{luc}(\MC(\R^d),L_1(\R^d))$, for all $t>0$,
$\xi_t(x) =\de \phi_t(Y)/\de Y(x)$ represents the unique solution of the equation
\[
\xi_t(x;z)=G_t(z,x) + \int_0^t ds \int G_{t-s} (z,y) \left( \int \frac{\de V_s(\phi_s,y)}{\de \phi_s(w)}\xi_s(x;w) \, dw
+ V_s(\phi_s,y)\xi_s(x;y)\right) dy
\]
\begin{equation}
\label{eq3thwellposMcVlasense}
- \int_0^t ds \int \left(\frac{\pa}{\pa y} G_{t-s} (z,y),
\int \frac{\de b_s(\phi_s,y)}{\de \phi_s(w)}\xi_s(x;w) \, dw
+ b_s(\phi_s,y)\xi_s(x;y)\right) dy.
\end{equation}
It satisfies  the Dirac initial condition $\xi_0=\de_x$ and has the bound
\begin{equation}
\label{eq4thwellposMcVlasense}
\|\xi_t(x,.)\|_{L_1(\R^d)}\le  E_{1/2}[C (2t+1) (2\la R(\la)+V+b)],
\end{equation}
with a constant $C$ depending on $m,M,T$.
Finally $\xi_t$ also solves the weak equation
 obtained by the formal differentiation of  (\ref{eq4thwellposfirstordernonl}):
\[
\frac{d}{dt}(f,\xi_t(x;.))=\left(\frac12 (A(.)\nabla, \nabla) f + (b_t(\phi_t, .),\nabla)f(.)+V_t(\phi_t, .)f(.), \xi_t(x;.)\right)
\]
\begin{equation}
\label{eq5thwellposMcVlasense}
+\int \int \frac{\de V_t(\phi_t,y)}{\de \phi_t(w)}\xi_t(x,w) f(y) \phi_t (y) \, dy dw
+\int \int \left( \frac{\de b_t(\phi_t,y)}{\de \phi_t(w)}\xi_t(x,w), \nabla f(y) \right) \phi_t(y)\, dy dw.
\end{equation}
\end{theorem}

\noindent{\it Proof.} 
The norm of the expression under the integral in (\ref{eq3thwellposMcVlasense}) is bounded by
\[
 (t-s)^{-1/2} C(m,M,T) [\la R(\la)+b+\sqrt t(\la R(\la)+V)].
\]

Furthermore, the norm of the difference of the expressions under the integral in (\ref{eq3thwellposMcVlasense}) for two
different functions $\phi^1$ and $\phi^2$ is bounded by
\[
 (t-s)^{-1/2} C(m,M,T) (\la L(\la)+R(\la))\|\phi^1-\phi^2\|_{C([0,s], L_1(\R^d))}(1+\sqrt t).
\]
Consequently the convergence of the iterations arising from equation (\ref{eq3thwellposMcVlasense})
again follows from the fixed point principle. Once the existence of bounded solution to (\ref{eq3thwellposMcVlasense})
is obtained, one shows that it represents the derivative of $\phi_t$ with respect to initial data
in the usual way, as is done for the equation with the Lipschitz continuous right hand side (see e.g. \cite{Kol1, Kol2}).

\hfill $\Box$  

The weak and mild representations (\ref{eq5thwellposMcVlasense}) and (\ref{eq3thwellposMcVlasense}) of the derivatives
 with respect to the initial conditions can be used to derive different kinds of regularity for these derivatives.
 In order to illustrate this claim, let us observe that the weak equation (\ref{eq5thwellposMcVlasense}) shows that the evolution
$\xi\mapsto \xi_t$ of the directional derivatives of the solutions $\phi_t$ is dual to the evolution in $C(\R^d)$
generated by the equation
\[
\dot f_t(z)=\frac12 (A(.)\nabla, \nabla ) f_t(z) + (b_t(\phi_t, z),\nabla)f(z)+V_t(\phi_t, z)f_t(z)
\]
\begin{equation}
\label{eq6thwellposMcVlasense}
+\int \int \frac{\de V_t(\phi_t,y)}{\de \phi_t(z)} f_t(y) \phi_t (y) \, dy
+\int \int \left( \frac{\de b_t(\phi_t,y)}{\de \phi_t(z)}, \nabla f_t(y) \right) \phi_t(y)\, dy.
\end{equation}

\begin{theorem}
\label{propdervarder}
 Under the assumptions of Theorem \ref{thwellposMcVlasense} let additionally
\[
b_t(\mu, .), V_t(\mu, .), \,\,  \frac{\de b_t(\mu,y)}{\de \mu(.)}, \frac{\de V_t(\mu,y)}{\de \mu(.)} \in C^1(\R^d)
\]
uniformly for all $y$ and $\mu \in \MC^+_{\le \la}(\R^d)$ with any $\la$, so that
\begin{equation}
\label{eq00propdervarder}
\|b_t(\mu, .)\|_{C^1(\R^d)}+\|V_t(\mu, .)\|_{C^1(\R^d)}
+\|\frac{\de b_t(\mu,y)}{\de \mu(.)}\|_{C^1(\R^d)}
+ \|\frac{\de V_t(\mu,y)}{\de \mu(.)}\|_{C^1(\R^d)}\le c_1(\la).
\end{equation}

Then the equation (\ref{eq6thwellposMcVlasense}) is well-posed in $C^2_{\infty}(\R^d)$
and it generates a backward propagator $\Phi^{t,s}$ acting strongly continuously in the space
$C_{\infty}(\R^d)$,  $C^1_{\infty}(\R^d)$,  $C^2_{\infty}(\R^d)$, so that
\[
\|\Phi^{t,s}\|_{\LC(C^1_{\infty}(\R^d))}\le e^{C(m,M)(s-t)} E_{1/2} [(V+b+2R(\la)) C(m,M) \sqrt{s-t}],
\]
\begin{equation}
\label{eq1propdervarder}
\|\Phi^{t,s}\|_{\LC(C^2_{\infty}(\R^d))}\le e^{C(m,M)(s-t)} E_{1/2} [c_1(\la) C(m,M) \sqrt{s-t}],
\end{equation}
with a constant $C(m,M)$ depending only on $m,M$.
Consequently, by duality, the weak equation (\ref{eq5thwellposMcVlasense}) generates a (forward) propagator
$(\Phi^{t,s})^*$ in $\MC(\R^d)$ that extends to bounded
propagators in the dual spaces $(C^1_{\infty}(\R^d))^*$ and  $(C^2_{\infty}(\R^d))^*$. In particular, the variational derivative $\xi_t(x;.)$
are twice differentiable in $x$, so that
\[
\pa \xi_t(x;.)/\pa x\in (C^1_{\infty}(\R^d))^*, \quad
\pa^2 \xi_t(x;.)/\pa x^2 \in (C^2_{\infty}(\R^d))^*,
\]
and
\[
\pa \xi_t(x;.)/\pa x|_{(C^1_{\infty}(\R^d))^*} \le e^{C(m,M)t} E_{1/2} [(V+b+2R(\la)) C(m,M) \sqrt t],
\]
\begin{equation}
\label{eq2propdervarder}
\pa \xi_t(x;.)/\pa x|_{(C^2_{\infty}(\R^d))^*} \le e^{C(m,M)t} E_{1/2} [c_1(\la) C(m,M) \sqrt t].
\end{equation}
\end{theorem}

\noindent{\it Proof.} 
 It is a consequence of the perturbation theory applied to equation (\ref{eq6thwellposMcVlasense}).
\hfill $\Box$  

Let us look at the second variational derivatives of the solutions of the McKean-Vlasov equations with respect to the initial data:
\[
\eta_t(x,z;.)=\de^2\phi_t(Y)/\de Y(x)\de Y(z).
\]

Differentiating equation  (\ref{eq5thwellposMcVlasense}) we obtain for $\eta$ the weak equation
\[
\frac{d}{dt}(f,\eta_t(x,z;.))
=\left(\frac12 (A(.)\nabla, \nabla) f + (b_t(\phi_t, .),\nabla)f(.)+V_t(\phi_t, .)f(.), \eta_t(x,z;.)\right)
\]
\[
+(f,q_t) +\int \int \frac{\de V_t(\phi_t,y)}{\de \phi_t(w)}\eta_t(x,z;w) f(y) \phi_t (y) \, dy dw
\]
\begin{equation}
\label{McVlassecder}
+\int \int \left( \frac{\de b_t(\phi_t,y)}{\de \phi_t(w)}\eta_t(x,z;w), \nabla f(y) \right) \phi_t(y)\, dy dw
\end{equation}
with $(f,q_t)$ being given by
\[
\int \int \left[\frac{\de V_t(\phi_t,y)}{\de \phi_t(w)} f(y)+\left( \frac{\de b_t(\phi_t,y)}{\de \phi_t(w)}, \nabla f(y) \right)\right]
[\xi_t(x;y)\xi_t(z;w)+ \xi_t(x;w)\xi_t(z;y)] \, dy dw
\]
\begin{equation}
\label{McVlassecder1}
\int \int \iint \left[\frac{\de^2 V_t(\phi_t,y)}{\de \phi_t(w)\de \phi_t(u)} f(y)
+\left( \frac{\de ^2 b_t(\phi_t,y)}{\de \phi_t(w)\de \phi_t(u)}, \nabla f(y) \right)\right]
\xi_t(x;w)\xi_t(z;u) \phi_t(y) \, dy dw du,
\end{equation}
which should be satisfied with the vanishing initial condition.

This is the same equation as (\ref{eq5thwellposMcVlasense}), but with the additional non-homogeneous term $(f, q_t)$. Thus its solution
can be expressed in terms of the propagators $(\Phi^{t,s})^*$ from Theorem \ref{propdervarder}:
\begin{equation}
\label{McVlassecder2}
\eta_t(x,z;.)=\int_0^t (\Phi^{0,s})^* q_s \, ds.
\end{equation}
Therefore, in spite of the formidably looking expression (\ref{McVlassecder1}), analysis of $\eta$ is more or less straightforward.

The structure of (\ref{McVlassecder1}) conveys an important message that for this analysis one needs the exotic spaces
$C^{2,k \times k} (\MC^+_{\le \la}(\R^d))$ introduced in the introduction.

From Theorem \ref{propdervarder} and formula
(\ref{McVlassecder2}) we can now derive the following consequence on the {\it second order sensitivity for the McKean-Vlasov diffusion}.
\index{sensitivity for McKean-Vlasov diffusion!second order}

\begin{theorem}
\label{propdervardersec}
(i) Under the assumptions of Theorem \ref{propdervarder} assume the existence of continuous
bounded second order variational derivatives:
 \begin{equation}
\label{eq0propdervardersec}
\sup_{y,w,u,t} |\frac{\de ^2 V_t(\mu,y)}{\de \mu(w)\de \mu(u)}| \le R_2(\la),
\quad
\sup_{y,w,u,t} \sum_j |\frac{\de ^2 b^j_t(\mu,y)}{\de \mu(w)\de \mu(u)}| \le R_2(\la).
\end{equation}
Then $\eta_t(x,z;.)$ is well-defined for any $t$ as an element of $(C^1(\R^d))^*$ and it has the following bound:
\begin{equation}
\label{eq1propdervardersec}
\|\eta_t(x,z;.)\|_{(C^1(\R^d))^*}
\le t  C(m,M,T) \la [R(\la)+R_2(\la)] \left(E_{1/2}[C(m,M,T) (2\la R(\la)+V+b)]\right)^3.
\end{equation}
with $\la=\|Y\|_{\MC(\R^d)}$.

(ii) Assuming additionally that
 \begin{equation}
\label{eq2propdervardersec}
\sup_{y,t} \|\frac{\de ^2 V_t(\mu,y)}{\de \mu(.)\de \mu(.)}\|_{C^{1\times 1}(\R^{2d})}\le R_3(\la),
\quad
\sup_{y,t} \sum_j\|\frac{\de ^2 b^j_t(\mu,y)}{\de \mu(.)\de \mu(.)}\|_{C^{1\times 1}(\R^{2d})}\le R_3(\la),
\end{equation}
 \begin{equation}
\label{eq3propdervardersec}
\sup_{w,t} \|\frac{\de V_t(\mu,.)}{\de \mu(w)}\|_{C^1(\R^d)} \le R_4(\la),
\quad
\sup_{w,t} \sum_j\|\frac{\de b^j_t(\mu,.)}{\de \mu(w)}\|_{C^1(\R^d)} \le R_4(\la),
\end{equation}
it follows that the derivatives of $\eta_t(x,z;.)$ with respect to $x$ and $z$ of order at most one are well-defined
 as elements of $(C^2(\R^d))^*$ and
 \[
\|\frac{\pa^{\al}}{\pa x^{\al}} \frac{\pa^{\be}}{\pa z^{\be}} \eta_t(x,z;.)\|_{(C^2(\R^d))^*}
\]
 \begin{equation}
\label{eq4propdervardersec}
\le  t  C(m,M,T) \la [R(\la)+R_2(\la)+ R_3(\la)+ R_4(\la)] \left(E_{1/2}[C(m,M,T) c_1(\la)]\right)^3
\end{equation}
for $\al,\be=0,1$.
\end{theorem}

\noindent{\it Proof.} 
(i) We have for $q_t$ in (\ref{McVlassecder2}) the estimate
\[
\|q_t\|_{(C^1(\R^d))^*} \le 4 \la [R(\la)+R_2(\la)] \|Y\|_{\MC(\R^d)} \|\xi_t(x;.)\|_{L_1(\R^d)}  \|\xi_t(z;.)\|_{L_1(\R^d)}
\]
\[
\le  4 \la [R(\la)+R_2(\la)]  \left(E_{1/2}[C(m,M,T)c_1(\la)]\right)^2,
\]
where (\ref{eq4thwellposMcVlasense}) was used. Hence (\ref{eq1propdervardersec}) follows by
the first estimate in (\ref{eq1propdervarder}).

(ii) The same as (i), but using the second estimate in (\ref{eq1propdervarder}).
\hfill $\Box$ 

\section{Sensitivity for Stochastic McKean-Vlasov Equations}

In this section we shall study the sensitivity of the following second
order stochastic  McKean-Vlasov equations

\begin{equation}
	\label{eqlimMacVlaSPDEsens}
	d(f, \mu_t)=(L_{t}(\mu_t)f, \mu_t) \, dt +(\si_{com}(.)\nabla f, \mu_t) \, dW_t,
\end{equation}
where
\begin{equation}
	\label{eqlimMacVlaSPDEsens1}
	L_{t}(\mu_t)f (x)= \frac12
	((\sigma_{com}\sigma_{com}^T+\sigma_{ind}\sigma_{ind}^T)(x)\nabla,\nabla)f (x)
	+(b(t,x, \mu_t),\nabla)f (x),
\end{equation}
$f \in C^2(\R^d)$, $W_t$ is a standard $1$-dimensional Brownian motion, functions $\si_{com}(x)$ and $\sigma_{ind}(x)$ are, respectively, $d \times 1$ vector-valued and $d \times d'$ matrix-valued continuous functions on $\R^d$, function $b(t,x,\mu_t)$ is  $d \times 1$ vector-valued continuous function on $[0,T] \times \R^d \times \MC(\R^d)$.

As we mentioned already it is more convenient to work with equation (\ref{eqlimMacVlaSPDEsens})
rewritten in the Stratonovich form as
\begin{equation}
	\label{eqlimMacVlaSPDEst2}
	d(f, \mu_t)=(L_t^{St}(\mu_t) f, \mu_t) \, dt +
	(\si_{com}(.) \nabla f, \mu_t) \, \circ dW_t,
\end{equation}
where
\[
L_t^{St}(\mu_t)f (x)=
\frac{1}{2}(\si_{ind}(x)\si_{ind}^T(x)\nabla,\nabla)f(x)
+\left((b(t,x,\mu_t)-\frac{1}{2}(\si_{com}(x)\nabla \si_{com}(x))),\nabla \right)f(x).
\]

The strong form of equation (\ref{eqlimMacVlaSPDEst2}) written for measures $\mu_t$ having densities $v_t$ is
\begin{equation}
	dv_t = {L_t^{St}}^{*}(v_t) v_t(x)\ dt - \nabla(\si_{com}(x)v_t(x)) \circ dW_t,
	\label{eqlimMacVlaSPDEst2str}
\end{equation}
where we denote the functions of measures and of their densities by the same letters, with some abuse of notations.
Here ${L_t^{St}}^{*}(v_t)$ is the dual operator to $L_t^{St}(v_t)$:
\begin{equation}
	\label{eqLstdual}
{L_t^{St}}^{*}(v_t)v_t (x)=
\frac{1}{2}(\nabla,\nabla (\si_{ind}(x)\si_{ind}^T(x)v_t(x)))
-(\nabla, ((b(t,x,v_t)-\frac{1}{2}(\si_{com}(x)\nabla \si_{com}(x)))v_t(x))).
\end{equation}

Let us introduce the following conditions:

(C1) $(\si_{ind}\si_{ind}^T)$ is uniformly elliptic matrix, that is
\[
m_{\si}^{-1}(\xi, \xi) \le (\si_{ind}(x)\si_{ind}^T(x)\xi, \xi) \le m_{\si} (\xi, \xi)
\]
for all $\xi, x\in \R^{d}$ and a constant $m_{\si}>0$;

(C2) Functions $\si_{ind}^{ij}(x), \si^i_{com}(x) \in C^4(\R^d)$;

(C3) Function $b(t,x,\mu)$ is continuous and bounded on $[0,T]\times \R^d\times \MC(\R^d)$, $b(t,.,\mu)\in C^2(\R^d)$, and $b$ is Lipshitz continuous as a function of $x$, uniformly in other variables;

(C4) The first and second order variational derivatives of $b(t,x,\mu)$ with respect to $\mu$ are well defined, bounded and
\[
b(t,x,.)\in (C^{2,1\times 1}\cap C^{1,2})(\MC_1(\R^d)).
\]

By using the method of stochastic characteristics, we are going to transfer equations (\ref{eqlimMacVlaSPDEst2}) and (\ref{eqlimMacVlaSPDEst2str}) to the non-stochastic equations with random coefficients thus making it possible to use the estimates obtained in the first part of this work.

Let us denote $\Omega f  = \si_{com}(x) \nabla f$.
The operator $\Omega' : v(x) \mapsto - \nabla(\si_{com}(x)v(x))$ is the dual operator to $\Omega$.
Let us denote
\[
B(x) = -\sum_{i=1}^{d}\frac{\pa \si_{com}^i}{\pa x^i}(\cdot).
\]

We can rewrite equation (\ref{eqlimMacVlaSPDEst2str}) as
\begin{equation}
\label{eqlimMacVlaSPDEst2strOm}
dv_t = {L_t^{St}}^{*}(v_t) v_t(x)\ dt + \Omega' v_t(x) \circ dW_t,
\end{equation}
where
$\Omega' v(x) = -(\si_{com}(x),\nabla)v(x)+B(x)v(x)$.

For $\si_{com}(x)\in C^2(\R^d)$, operator $\Omega'$ generates the contraction
group $e^{t\Om'}$ in $C(\R^d)$, so that $e^{t\Om'}v_0(x)$ is the unique solution to the equation
\[
\frac{\pa v}{\pa t}=\Om' v
\]
with the initial condition $v(0,x)=v_0(x)$. Explicitly,
\begin{equation}
	\label{eqmethcharac}
	e^{t\Om'}v_0(x)=v_0(Z(t,x))G(t,x), \quad t \in \R,
\end{equation}
where $Z(t,x)$ is the unique solution to the problem
\begin{equation}\label{ode Y}
	\frac{\partial Z}{\partial t}(t, x) = -\si_{com}(Z(t,x)), \qquad Z(0, x) = x,
\end{equation}
 and
\[
G(t,x)=\exp \{ \int_0^t B(Z(s,x)) \, ds \}.
\]
In particular, $G$ has the properties:
\[
G(-t,x)=G^{-1} (t, Z(-t,x))=\exp \{ \int_0^{-t} B(Z(s,x)) \, ds \}=\exp \{ -\int_{-t}^0 B(Z(s,x)) \, ds \},
\]

Thus from the standard formulas for the derivatives of ODEs with respect to initial data we get the following.

\begin{prop}\label{propY}
Under condition (C2) the function $Z(t, \cdot)$ is a flow of $C^3$-diffeomorphisms on $\R^d$ such that all derivatives of $Z(t,x)$ of order up to and including $3$
with respect to both its variables are uniformly bounded by $C\exp\{Ct\}$ with some constant $C$.
\end{prop}

Making in (\ref{eqlimMacVlaSPDEst2strOm}) the change of function $v$ to $g=\exp \{-\Om' W_t\} v$ leads to the equation
\begin{equation}
\label{eqstochMcKeanVltran}
\dot g_t =\exp \{-\Om' W_t\} {L_{t,\exp \{\Om' W_t\} g_t}^{{St}^*}}\exp \{\Om' W_t\} g_t,
\end{equation}
or in the weak form
\begin{equation}
\label{eqstochMcKeanVlweaktran}
\frac{d}{dt} (f,g_t) = (\tilde L_{t,\exp \{\Om' W_t\}g_t} f,g)=(\exp \{\Om W_t\} L_{t,\exp \{\Om' W_t\}g_t}^{St} \exp \{-\Om W_t\}f, g).
\end{equation}

We have
\begin{equation}
\label{exp Om Wt_g}
\exp\{\Om' W_t\}p(x) = G(W_t,x) p(Z(W_t,x)), \quad \exp\{\Om W_t\}p(x)= p(Z(-W_t,x)).
\end{equation}

To calculate $\tilde L$ it is convenient to first calculate the general effect of 'dressing' leading to the following result.
\begin{lemma}
\label{lemmadress}
Let
\[
\Lambda p(x)=\frac12 (\al(x)\nabla, \nabla) p(x)+\be (x) \nabla p(x).
\]
Then
\[
\exp \{\Om W_t\} \Lambda \exp \{-\Om W_t\}p(x)
\]
\[
=\frac12
(\frac{\pa Z(W_t,z)}{\pa z}|_{z=Z(-W_t,x)})^T \al(Z(-W_t,x))\frac{\pa Z(W_t,z)}{\pa z}|_{z=Z(-W_t,x)}
\nabla, \nabla) p(x)
\]
\[
+\left(\frac{\pa Z(W_t,z)}{\pa z}|_{z=Z(-W_t,x)}  \be (Z(-W_t,x)), \nabla p(x)\right)
\]
\[
+\left(\frac{\pa^2 Z(W_t,z)}{\pa z^2}|_{z=Z(-W_t,x)} \al(Z(-W_t,x)), \nabla p(x)\right).
\]
\end{lemma}

Applying this lemma yields
\begin{equation}
\tilde L_{t,\exp \{\Om' W_t\}g_t}f(x)
=\frac{1}{2}(\tilde A \nabla,\nabla)f(x)
+\left(\tilde{b}(t,x,[g]),\nabla\right)f(x),
\end{equation}
where
\begin{eqnarray}
\label{sitilda}
\tilde A^{ij}(x) = \left.\sum_{k,l=1}^{d}\left(\si_{ind}\si_{ind}^T\right)^{k l} \!(z)\frac{\pa Z_i}{\pa z^k}(W_t,z) \frac{\pa Z_j}{\pa z^l} (W_t,z) \right|_{z=Z(-W_t,x)},
\end{eqnarray}
\begin{eqnarray*}
& \tilde{b}^i(t,x,[g]) = \left.\displaystyle\sum_{k=1}^{d}\left( b^k(t,z,[\exp \{\Om' W_t\}g_t])
-\frac{1}{2} \sum_{l=1}^{d} (\si_{com}^l\frac{\pa \si_{com}^k}{\pa x^l})(z)\right)\frac{\pa Z_i}{\pa z^k}(W_t,z) \right|_{z=Z(-W_t,x)}& \\  %
&\displaystyle+\frac{1}{2} \sum_{k,l=1}^{d}\left(\si_{ind}\si_{ind}^T\right)^{k l} \!(z) \left.\frac{\pa^2 Z_i}{\pa z^k \pa z^l}(W_t,z)\right|_{z=Z(-W_t,x)}. &
\end{eqnarray*}

Let us denote ${\tilde L}'_{t,\exp \{\Om' W_t\}g_t}$ the dual operator of the operator ${\tilde L}_{t,\exp \{\Om' W_t\}g_t}$.

Then we can rewrite equation (\ref{eqstochMcKeanVltran}) as
\begin{equation}
\label{eqstochMcKeanVltranstrong}
\dot g_t={\tilde L}'_{t,\exp \{\Om' W_t\}g_t}g_t
=\frac{1}{2}(\nabla,\nabla(\tilde{\si}(x)g_t(x)))
-\left(\nabla,(\tilde{b}(t,x,[g])g_t(x))\right).
\end{equation}

The ellipticity is the crucial property of a diffusion operator. By Proposition \ref{propY},
\begin{equation}
\label{eqderbo}
\|\tilde{\si}(x)\|\le Ce^{C|W_t|}, \quad \|\tilde{\si}^{-1}(x)\|\le Ce^{C|W_t|}
\end{equation}
for all $x$. Uniform boundedness is valid only in some particular cases.

\begin{lemma}
\label{lemmaboundederinit}
Under condition (C2) if either (i) $\si_{com}(x) =(\si_{com}^1(x_1), \cdots , \si_{com}^d(x_d))$
or (ii) all solutions to the equation $\dot Z=\si_{com}(Z)$ are periodic with a fixed period
$T_{per}$, then the derivatives $\pa^k Z(t,x)/\pa x^k$, $k=1,2,3$, are uniformly bounded in $t$ and $x$, and consequently the coefficients
$\tilde A$ and $\tilde b$ have the same uniform bounds as the coefficients $\si_{ind}\si_{ind}^T$ and $b$ of the 'undressed equation'.
\end{lemma}

\noindent{\it Proof.} 
Under (i) the system $\dot Z=\si_{com}(Z)$ decomposes into one-dimensional equations, where the required boundedness is easy to show.
Under (ii) $Z(t+nT_{per}, x)= Z(t,x)$ and hence the bounds for the derivatives for arbitrary times are reduced to the bounds for $T<T_{per}$,
where they follow from Proposition \ref{propY}.
\hfill $\Box$ 

Equation (\ref{eqstochMcKeanVltranstrong}) is a particular case of equation (\ref{eqdifnonlgen11}), so that the theory of the previous section applies.

\begin{theorem}
	\label{thwellposMcVlaSt}
(i) Under assumptions (C1)-(C4) the mild form of the Cauchy problem for equation (\ref{eqstochMcKeanVltranstrong})
with initial condition $Y\in \MC^+(\R^d)$ is well-posed
almost surely, that is, the equation
\begin{equation}
		\label{eq1athwellposMcVla2}
		g_t(x)=\int G_t(x,y) Y(y)(dy)	- \int_0^t ds \int \left(\frac{\pa}{\pa y} G_{(t-s)} (x,y), \tilde b(s,y, [g_s]) g_s (y) \right) dy
	\end{equation}
has the unique bounded nonnegative solution $g_. \in C((0,T], L_1(\R^d))$
 such that $g_t \to Y$ weakly, as $t \to 0$ (and $g_t\to g_0$ in $L^1(\R^d)$ if $Y$ has a density $g_0$), and $\|g_t\|\le \|Y\|$. This implies
 the almost sure well-posedness  for equation (\ref{eqlimMacVlaSPDEst2str}) with the same bound $\|v_t\|\le \|Y\|$.
Solutions $g_t$ to the mild equation (\ref{eq1athwellposMcVla2}) also solve the corresponding weak equation
	\begin{equation}
		\label{eq4thwellposfirstordernonl2}
		\frac{d}{dt}(f,g_t)=\left(\frac12 (\tilde A \nabla, \nabla) f + (\tilde b(t,x, [g]),\nabla)f(x), g_t\right),
		\quad g_0=V\in \MC^+_{\le \la}(\R^d).
	\end{equation}

(ii)	If additionally assumptions of Lemma \ref{lemmaboundederinit} hold, then
for any two solutions $g_t^1$ and $g_t^2$ of (\ref{eq1athwellposMcVla2}) and the corresponding solutions of (\ref{eqlimMacVlaSPDEst2str})
	with the initial conditions $Y^1,Y^2$ the estimate
	\begin{equation}
		\label{eq2thwellposMcVlaSt}
		\|g_t^1 -g_t^2\|_{L_1(R^d)} \le \|Y^1-Y^2\|_{\MC(R^d)}E_{1/2} (\ka(T) \sqrt t),
	\end{equation}
\begin{equation}
		\label{eq2thwellposMcVlaa}
		\|v_t^1 -v_t^2\|_{L_1(R^d)} \le \|Y^1-Y^2\|_{\MC(R^d)}E_{1/2} (\ka(T) \sqrt t),
	\end{equation}
	holds, where
	\[
	\ka(T) =  \bar \ka (\sqrt T+1)\|Y\|,
	\]
and $\bar \ka$ depends on the bounds of the derivatives in conditions (C1)-C4).
\end{theorem}

\noindent{\it Proof.} 
The results for $g_t$ are the consequences of Theorem \ref{thwellposMcVla}.
For $v_t$ they follow from the formula $v_t=\exp\{\Om'W_t\} g_t$.
\hfill $\Box$ 

For application of this results it is impostant to observe the following consequence of formula (\ref{eq2thwellposMcVlaSt}):

	\begin{equation}
	\label{eq2thwellposMcVla2}
	\E \|g_t^1 -g_t^2\|_{L_1(R^d)} \le \|g_0^1-g_0^2\|_{\MC(R^d)} \exp (CT)
	\end{equation}
with some constant $C$.

Again directly applying Theorem \ref{thwellposMcVlasense} yields the following sensitivity result for the
stochastic McKean-Vlasov equation (\ref{eqstochMcKeanVltranstrong}).

\begin{theorem}
\label{thwellposMcVlasenseSt}
Under the Conditions (C1)-(C4) and the assumptions of Lemma \ref{lemmaboundederinit},
all results and estimates of Theorem \ref{thwellposMcVlasense} (with $V=0$) hold for equation (\ref{eq1athwellposMcVla2})
and hence for equation (\ref{eqstochMcKeanVltranstrong}).
\end{theorem}

Similarly Theorem \ref{propdervarder} on the regularity is automatically transferred to equation
(\ref{eqstochMcKeanVltranstrong}). Assumptions (C1)-(C4) and the assumptions of Lemma \ref{lemmaboundederinit}
ensure that the conditions of Theorem \ref{propdervarder} hold for equation (\ref{eq1athwellposMcVla2}).

Similarly by Theorem \ref{propdervardersec} we can now derive the following result on the second order sensitivity for the McKean-Vlasov equation (\ref{eqstochMcKeanVltranstrong}).

Denote by
\[
\eta_t(x,z;.)=\de^2g_t(Y)/\de Y(x)\de Y(z), \quad \tilde \eta_t(x,z;.)=\de^2v_t(Y)/\de Y(x)\de Y(z)
\]
the second variational derivatives of the solutions of the McKean-Vlasov SDEs with respect to the initial data.

\begin{theorem}
	\label{propdervardersecSt}
Under Conditions (C1)- (C4), $\eta_t(x,z;.)$ is well-defined for any $t$ as an element of $(C^1(\R^d))^*$ and
	\begin{equation}
	\label{eq1propdervardersecSt}
	\|\eta_t(x,z;.)\|_{(C^1(\R^d))^*}
	\le t  C(m,M,T) \la [R(\la)+R_2(\la)] \left(E_{1/2}[C(m,M,T) (\la R(\la)+b_1)]\right)^3.
	\end{equation}
	with $\la=\|Y\|_{\MC(\R^d)}$. 	
Moreover, the derivatives of $\eta_t(x,z;.)$ with respect to $x$ and $z$ of order at most one are well-defined
	as elements of $(C^2(\R^d))^*$ and
	\[
	\|\frac{\pa^{\al}}{\pa x^{\al}} \frac{\pa^{\be}}{\pa z^{\be}} \eta_t(x,z;.)\|_{(C^2(\R^d))^*}
	\]
	\begin{equation}
	\label{eq4propdervardersecSt}
	\le  t  C(m,M,T) \la [R(\la)+R_2(\la)+ R_3(\la)+ R_4(\la)] \left(E_{1/2}[C(m,M,T) c_1(\la)]\right)^3
	\end{equation}
	for $\al,\be=0,1$.
The same estimates hold for $\tilde \eta_t$.
\end{theorem}

\subsection{Multidimensional Common Noise}

In this pharagraph we shall study the well-posedness of the stochastic  McKean-Vlasov equations with multidimensional common noise.
Namely, the following McKean-Vlasov SPDE
\begin{equation}
\label{eqlimMacVlaSPDEmulti}
d(f, \mu_t)=(L_{t}(\mu_t)f, \mu_t) \, dt +(\si_{com}\nabla f, \mu_t) \, dW_t,
\end{equation}
where
\begin{equation}
\label{eqlimMacVlaSPDElmulti}
L_{t}(\mu_t)f (x)= \frac12
((\sigma_{com}\sigma_{com}^T+\sigma_{ind}\sigma_{ind}^T(x))\nabla,\nabla)f (x)
+(b(t,x, \mu_t),\nabla)f (x),
\end{equation}
$f \in C^2(\R^d)$, $W_t$ is a standard $d''$-dimensional Brownian motion, $\sigma_{ind}(x)$ is a $d \times d'$ matrix-valued continuous function on $\R^d$, function $b(t,x,\mu_t)$ is  $d \times 1$ vector-valued continuous function on $[0,T] \times \R^d \times \MC(\R^d)$ satisfying Conditions (C1)-(C4). As a crucial simplifying assumption as compared with previous study we assume that
$\si_{com}$ is a constant $d\times d''$ matrix.

We can rewrite the equation (\ref{eqlimMacVlaSPDEmulti}) in the Stratonovich form as

\begin{equation}
\label{eqlimMacVlaSPDEst2multi}
d(f, \mu_t)=(L_t^{St}(\mu_t) f, \mu_t) \, dt +
(\si_{com} \nabla f, \mu_t) \, \circ dW_t,
\end{equation}
where
\[
L_t^{St}(\mu_t)f (x)=
\frac{1}{2}(\si_{ind}(x)\si_{ind}^T(x)\nabla,\nabla)f(x)
+\left(b(t,x,\mu_t),\nabla \right)f(x).
\]

The strong form of equation (\ref{eqlimMacVlaSPDEst2multi}) written for measures $\mu_t$ having densities $v_t$ is
\begin{equation}
dv_t = {L_t^{St}}^{*}(v_t) v_t(x)\ dt - (\si_{com}, \nabla)v_t(x) \circ dW_t.
\label{eqlimMacVlaSPDEst2strmulti}
\end{equation}
where
\begin{equation}
\label{eqLstdualm}
{L_t^{St}}^{*}(v_t)v_t (x)=
\frac{1}{2}(\nabla,\nabla (\si_{ind}(x)\si_{ind}^T(x)v_t(x)))
-(\nabla, b(t,x,v_t)v_t(x)).
\end{equation}

As above for one-dimensional noise, let us denote
\[
\Omega v(x) = (\si_{com},\nabla)v(x)=\{\Omega^j v(x)\}= \{ \sum_k \si_{com}^{jk} \nabla_k v(x)\}.
\]
However now, this is a vector-valued operator generating $d''$ semigroups
\[
T_k(t)v(x)=v(x+\si_{com}^{k .}t)=v(x_j+\si_{com}^{kj}t),
\]
solving the Cauchy problems for the equations
\[
\frac{\pa v}{\pa t}(t,x) =\sum_j \si_{com}^{jk} \frac{\pa v}{\pa x_j}(t,x).
\]
These semigroups commute and define the action $T(t_1, \cdots , t_{d''})$ of $\R^{d''}$ on $(\R^d)$ by the formula
\[
T(t_1, \cdots , t_{d''}): v(x) \mapsto v(x+ \sum_k \si_{com}^{k.} t_k).
\]
Moreover,
the operators
\begin{equation}
\label{ivpsolmult}
T(W_t) v_0(x)= v_t(x) = v_0(x+\si_{com} W_t)=v_0(x+ \sum_k \si_{com}^{k.} W_k(t))
\end{equation}
solve the initial value problem for the first order SPDE
\begin{equation}
\label{ivpomegam}
dv_t = (\si_{com},\nabla) v_t(x) \circ dW_t.
\end{equation}

The operator
\[
\Omega' v(x) = -\Omega v(x)=(-\si_{com},\nabla)v(x),
\]
defines the flow that differs from $T(t_1, \cdots , t_{d''})$ by the inversion of time.

We can rewrite the equation (\ref{eqlimMacVlaSPDEst2strmulti}) as
\begin{equation}
\label{eqlimMacVlaSPDEst2strOmmulti}
dv_t = {L_t^{St}}^{*}(v_t) v_t(x)\ dt + \Omega' v_t(x) \circ dW_t.
\end{equation}

Making in (\ref{eqlimMacVlaSPDEst2strOmmulti}) the change of function $v$ to $g=T(W_t)v=v(x+\si_{com} W_t)$ leads to the equation
\begin{equation}
	\label{eqstochMcKeanVltranmulti}
	\dot g_t = \left[{L_{t, g_t(x-\si_{com} W_t)}^{{St}^*}} g_t(x-\si_{com} W_t)\right](x+\si_{com} W_t),
\end{equation}
or in the weak form
\begin{equation}
	\label{eqstochMcKeanVlweaktranmulti}
	\frac{d}{dt} (f,g_t) = (\tilde L_{t,g_t(x-\si_{com} W_t)} f,g)=([L_{t,g_t(x-\si_{com} W_t)}^{St} f(x+\si_{com} W_t)](x-\si_{com} W_t), g).
\end{equation}

Then we can prove by a direct computation that $\tilde L_{t,g_t(x-\si_{com} W_t)}$ in (\ref{eqstochMcKeanVlweaktranmulti}) is a second order differential operator represented by
\begin{equation}
\tilde L_{t,g_t(x-\si_{com} W_t)}f(x)
=\frac{1}{2}(\tilde A \nabla,\nabla)f(x)
+\left(\tilde{b}(t,x,[g]),\nabla\right)f(x),
\end{equation}
where
\[
\tilde A^{ij}(x) =(\si_{ind}\si_{ind}^T)(x +\si_{com} W_t),
\]
\[
\tilde{b}(t,x,[g]) =b(t,(x+\si_{com} W_t),[g_t(x-\si_{com} W_t)]).
\]

Let us denote ${\tilde L}'_{t,g_t(x-\si_{com} W_t}$ the dual operator of the operator ${\tilde L}_{t,g_t(x-\si_{com} W_t}$.

Then we can rewrite the equation (\ref{eqstochMcKeanVltranmulti}) as follows

\begin{equation}
\label{eqstochMcKeanVltranstrongmulti}
\dot g_t={\tilde L}'_{t,g_t(x-\si_{com} W_t)} g_t
=\frac{1}{2}(\nabla,\nabla(\tilde A(x)g_t(x)))
-\left(\nabla,(\tilde{b}(t,x,[g])g_t(x))\right).
\end{equation}

Equation (\ref{eqstochMcKeanVltranstrongmulti}) is a particular case of the equation (\ref{eqstochMcKeanVltranstrong}), so that
the sensitivity results given for the case of $1-$dimensional $W_t$ have direct extensions to the present case of multidimensional $W_t$ with some simplifications.

\section{Acknowledgments}
The research has been partially supported by the Ministry of Education and Science of the Russian Federation (Grant No. 1.6069.2017/8.9)




\end{document}